\documentclass{amsart}

\newtheorem{theorem}{Theorem}[section]
\newtheorem{lemma}[theorem]{Lemma}
\newtheorem{proposition}[theorem]{Proposition}
\newtheorem{corollary}[theorem]{Corollary}

\theoremstyle{definition}

\newtheorem{example}[theorem]{Example}
\newtheorem{examples}[theorem]{Examples}
\newtheorem{remark}[theorem]{Remark}

\theoremstyle{remark}

\def\boxtext#1{%
\vbox{%
\hrule
\hbox{\strut \vrule{} #1 \vrule}%
\hrule
}%
}

\usepackage{amscd,amssymb}

\begin{document}

\title[Constants of derivatives and Taylor expansions]
{Constants of formal derivatives\\ of non-associative algebras,\\
Taylor expansions and applications}

\author[Vesselin Drensky and Ralf Holtkamp]
{Vesselin Drensky and Ralf Holtkamp}
\address{Institute of Mathematics and Informatics,
Bulgarian Academy of Sciences,
          1113 Sofia, Bulgaria}
\email{drensky@math.bas.bg}
\address{Fakult\"at f\"ur Mathematik, Ruhr-Universit\"at,
          44780 Bochum, Germany}

\curraddr{} \email{ralf.holtkamp@ruhr-uni-bochum.de}

\thanks
{The work of the first author was partially supported by Grant
MM-1106/2001 of
the Bulgarian Foundation for Scientific Research.}


\begin{abstract} We study unitary multigraded non-associative algebras
$R$ generated by an ordered set $X$ over a field $K$ of
characteristic 0 such that the mappings $\partial_k:
x_l\to\delta_{kl}$, $x_k,x_l\in X$, can be extended to derivations
of $R$. The class of these algebras is quite large and includes
free associative and Jordan algebras, absolutely free
(non-associative) algebras, relatively free algebras in varieties
of algebras, universal enveloping algebras of multigraded Lie
algebras, etc.

There are Taylor-like formulas for $R$: Each element of $R$ can be
uniquely presented as a sum of elements of the form
$(\cdots(r_0x_{j_1})\cdots x_{j_{n-1}}) x_{j_n}$, where $r_0$ is a
constant (i.~e. $\partial_k(r_0)=0$ for all $x_k\in X$) and
$j_1\leq\cdots\leq j_{n-1}\leq j_n$. We present methods for the
description of the algebra of constants, including an approach via
representation theory of the general linear or symmetric groups.
 As an application of the Taylor expansion for non-associative algebras,
 we consider the solutions of ordinary
linear differential equations with constant coefficients from the
base field.

\end{abstract}

\maketitle

\section*{Introduction}

Let $K\{X\}$ be the (absolutely) free unitary non-associative algebra freely
generated by an ordered set $X$ over a field $K$ of characteristic
0 \ ($K\{X\}$ is also called the free magma algebra over $X$).
We consider factor algebras $R=K\{X\}/I$, where the ideal $I$
is multigraded and invariant under the formal derivatives
$\partial_k=\frac{\partial}{\partial x_k}$, $x_k\in X$. We shall
denote the generators $x_j+I$ of $R$, $x_j\in X$, with the same
symbols $x_j$. Since $\partial_k(I)\subseteq I$, the action of
$\partial_k$ induces a derivation of $R$. Let $R_0$ be the algebra
of constants of $R$, i.~e. the set of all $r_0\in R$ such that
$\partial_k(r_0)=0$ for all $x_k\in X$.
\par
The description of $R_0$ is known in many cases. For example, an old
result of
Falk \cite{F} describes the case when $R=K\langle X\rangle$ is the
free associative algebra. Then $R_0$ is generated by all Lie
commutators
$[\ldots[x_{j_1},x_{j_2}],\ldots],x_{j_n}]$, $n\geq 2$. Specht
\cite{Sp}
applied products of such commutators in the study of algebras with
polynomial
identities, see the book of one of the authors \cite{D3} for further
application to the theory of PI-algebras. It is known, \cite{G3},
that in this case the algebra $R_0$ is free,
see also \cite{DK} for an explicit basis of $R_0$.
\par
Any element $r\in R$ can be expressed in unique way in the form $$
r=\sum (\cdots (r_jx_{j_1})\cdots x_{j_{n-1}})x_{j_n},\ r_j\in
R_0,\ j_1\leq\cdots\leq j_{n-1}\leq j_n. $$

Such Taylor-like formulas were used for relatively free
associative and Jordan algebras, see \cite{D1, D2}, and on skew
polynomial constructions for the free associative algebra, see
\cite{G1}, and on free magma algebras $K\{X\}$, see
\cite{GH}. We obtain more general results about factor algebras
$R=K\{X\}/I$ in Section 2.
\par
In Section 3, we consider the algebra of constants in the case of
$K\{X\}$ and the case of free commutative non-associative
algebras. We give different proofs for results of \cite{GH}, and
we sketch methods to obtain a concrete description, including a
method using representation theory of the general linear (or
symmetric) groups, which follows ideas of Regev \cite{R} and one
of the authors \cite{D1}.

In the last section we consider ordinary linear differential
equations of the form $$
y^{(n)}+a_1y^{(n-1)}+\cdots+a_{n-1}y'+a_ny=f(x), $$ where the
coefficients $a_1,\ldots,a_{n-1},a_n$ are constants from the field
$K$ and $f(x)$ is a formal power series from the completion $\hat
R$ of the algebra $R=K\{x\}/I$ where $I$ is a homogeneous ideal
invariant under the formal derivative $d_x=d/dx$. In this setup,
we establish an analogue of the well know result about the general
form of the solutions of such equations in terms of real and
complex functions. A special case of our considerations is the
non-associative exponential function defined in \cite{DG}.

\section{Conventions}

We fix a field $K$ of characteristic 0 and an ordered set
$X=\{x_j\mid j\in J\}$.

Let $K\{X\}$ be the absolutely free non-associative and
noncommutative unitary $K$-algebra freely generated by the set
$X$. In the
following, by an algebra we just mean a $K$-vector space $V$
together with a (not necessarily associative) binary operation
$\cdot: V\times V\to V$. Every algebra is a factor algebra
$R=K\{X\}/I$, where $I$ is an ideal of $K\{X\}$.

 As a vector space $K\{X\}$ has a
basis consisting of all non-associative words in the alphabet $X$.
For example, we make a difference between $(x_1x_2)x_3$ and
$x_1(x_2x_3)$ and even between $(xx)x$ and $x(xx)$. Words of
length $n$ correspond to planar binary trees with $n$ leaves, see
e.~g. \cite{H}, \cite{GH}.

We call the elements of $K\{X\}$ polynomials in non-associative
and noncommutative variables. We omit the parentheses when the
products are left normed. For example, $uvw=(uv)w$ and
$x^n=(x^{n-1})x$. The algebra $K\{X\}$ has a natural multigrading,
counting the degree $\text{deg}_j(u)$ of any non-associative word
$u$ with respect to each free generator $x_j\in X$. Every mapping
$\delta:X\to K\{X\}$ can be extended to a derivation of $K\{X\}$:
If $u,v$ are monomials in $K\{X\}$, then we define $\delta(uv)$
inductively by $\delta(uv)=\delta(u)v+u\delta(v)$.
\par
In this paper we are interested in the formal partial derivatives
$\partial_k$ defined by
\[
\partial_k(x_l):=\frac{\partial x_l}{\partial x_k}:=\delta_{kl}=
\begin{cases}1&, \text{ if } k=l,\\
0&, \text{ if } k\not=l.
\end{cases}
\]
If $I$ is any multihomogeneous ideal of $K\{X\}$ such that
$\partial_k(I)\subseteq I$ for all $x_k\in X$, then the factor
algebra $R=K\{X\}/I$ admits the formal partial derivatives. In the
sequel we shall consider such ideals $I$ and such algebras $R$
only. We shall use the same symbols $x_j,\partial_k$, etc. for the
images in $R$ of the corresponding objects in $K\{X\}$. We denote
by $R_0$ the subalgebra of $R$ consisting of all constants, i.~e.
$R_0$ is the intersection of the kernels $\text{Ker }\partial_k$,
$x_k\in X$.

Let $u\in R$. By
 $\lambda_u$ and $\rho_u$, respectively, we denote the operators
 of left and right multiplication defined by
 $$
\lambda_u:v\to u v,\ \rho_u:v\to v u,\ v\in R. $$ The left and
right multiplications $\lambda_u$ and $\rho_v$ (all $u,v$)
generate a subalgebra ${\mathcal M}(R)$ of the (associative)
algebra of linear operators on the vector space $R$. The
associative algebra ${\mathcal M}(R)$ is called the
 algebra of multiplications of $R$.
The algebra ${\mathcal M}(R)$ inherits the multigrading of $R$ and
for any non-associative word $u\in R$ we have
$\text{deg}_j(\lambda_u)=\text{deg}_j(\rho_u)=\text{deg}_j(u)$,
$x_j\in X$.

\section{Taylor Formulas and Constants of Derivatives}

The following statement should be considered as a non-associative
and noncommutative analogue of the Taylor formula.

\begin{proposition} \label{Taylor}
Let $I$ be an ideal of $K\{X\}$ which is invariant
under the formal partial derivatives. Let $R_0$ be the subalgebra of
all constants in the factor algebra $R=K\{X\}/I$,
i.~e. $r_0\in R_0$ if and only if $\partial r_0/\partial x_k=0$
for all $x_k\in X$. Then every element $r\in R$ can be expressed
in a unique way in the form
\[
r=\sum (\cdots(r_jx_{j_1})\cdots x_{j_{n-1}})x_{j_n},
\]
where $r_j\in R_0$, $j=(j_1,\ldots,j_n)$, are constants and
$x_{j_1}\leq \cdots\leq x_{j_{n-1}}\leq x_{j_n}$.

If $r$ depends on the variables $x_1<\cdots<x_m$ only, then
\[
r=\sum r_a\rho_1^{a_1}\cdots \rho_m^{a_m},\ r_a\in R_0.
\]

\end{proposition}

\begin{proof}
For $R=K\{X\}$, this is  Proposition 3.1 of \cite{GH}. The proof
for the case $R=K\{X\}/I$ is similar. For the sake of
completeness, we sketch it, following the main steps of the proof
of \cite{D2} Proposition 1.5.

 Let $r\in R$ depend on
the variables $x_1<\cdots<x_m$ only. Consider the element
\[
r_0=r-\frac{\partial r}{\partial x_m}\frac{\rho_m}{1!}+
\frac{\partial^2 r}{\partial x_m^2}\frac{\rho_m^2}{2!}-
\frac{\partial^3 r}{\partial x_m^3}\frac{\rho_m^3}{3!}+\cdots
\]
It is easy to see that $\partial r_0/\partial x_m=0$. Since the total
degree of
$\partial^pr/\partial x_m^p$ is lower than the degree of $r$,
by induction we obtain that the derivatives $\partial^pr/\partial
x_m^p$
already have the desired expression
\[
\frac{\partial^pr}{\partial x_m^p}=\sum r_{pb}\rho_m^b,
\]
with $\partial r_{pb}/\partial x_m=0$.
Hence $r$ has the form
\[
r=\sum_{a\geq 0}r_a\rho^a_m,\
\frac{\partial r_a}{\partial x_m}=0.
\]
Continuing with the next variables $x_{m-1},x_{m-2},\ldots,x_1$, we
obtain the presentation.
\par
In order to see that the presentation is unique, let
\[
r=\sum r_a\rho_1^{a_1}\cdots \rho_m^{a_m}=0,\
r_a\in R_0,
\]
with some $r_a\not=0$. We choose the maximal $m$-tuple
$(k_1,\ldots,k_m)$ among all indices $(a_1,\ldots,a_m)$ with
$r_a\not=0$, with respect to the lexicographical order. Then
\[
0=\frac{\partial^{k_1+\cdots+k_m}r}{\partial x_1^{k_1}\cdots\partial
x_m^{k_m}}
=r_k\not=0,
\]
which is a contradiction.

\end{proof}

\begin{corollary}\label{ideals}
A multihomogeneous ideal $I$ of $K\{X\}$ is invariant under all formal
partial derivatives if and only if $I$ can be generated by constants.
\end{corollary}

\begin{proof}
If $I$ is generated by constants $r_j$, $j\in J$, then any element of
$I$
has the form
\[
r=\sum \alpha_{ju}r_j\mu_{u_1}\cdots\mu_{u_n},
\]
where $\alpha_{ju}\in K$ and $\mu_{u_p}$ is the operator of
left or right multiplication by the monomial $u_p$. Since
\[
\frac{\partial r}{\partial x_k}=\sum\alpha_{ju}\left(
\frac{\partial r_j}{\partial x_k}\mu_{u_1}\cdots\mu_{u_n}
+\sum_{p=1}^nr_j\mu_{u_1}\cdots\mu_{\partial u_p/\partial
x_k}\cdots\mu_{u_n}\right)
\]
and $\partial r_j/\partial x_k=0$, we obtain that $\partial r/\partial
x_k$
also belongs to $I$.
\par
In the opposite direction, let, in the notation of Proposition
\ref{Taylor},
\[
r=\sum r_a\rho_1^{a_1}\cdots \rho_m^{a_m}=0,\
r_a\in R_0,
\]
be a generator of the ideal $I$ invariant under all partial
derivatives. We choose the maximal $m$-tuple $(k_1,\ldots,k_m)$
among all indices $(a_1,\ldots,a_m)$ with $r_a\not=0$, with
respect to the lexicographical order, and obtain
\[
\frac{\partial^{k_1+\cdots+k_m}r}{\partial x_1^{k_1}\cdots\partial
x_m^{k_m}}
=r_k\in I.
\]
Hence we can replace the generator $r$ of $I$ with the constants
$r_a$ which completes the proof.
\end{proof}

\begin{remark}
The same arguments as above give that Corollary \ref{ideals}
is true if we replace $K\{X\}$ with a factor algebra $R=K\{X\}/I$
modulo a multihomogeneous ideal $I$ of $K\{X\}$
which is invariant under the formal derivatives.
\end{remark}

Recall, that if $V=\bigoplus V^{(n_1,\ldots,n_m)}$ is a
multigraded vector space where $V^{(n_1,\ldots,n_m)}$ is the
multihomogeneous component of degree $(n_1,\ldots,n_m)$, then the
Hilbert (or Poincar\'e) series of $V$ is defined by
\[
\text{\rm Hilb}(V,t_1,\ldots,t_m) =\sum
\text{dim }V^{(n_1,\ldots,n_m)}t_1^{n_1}\cdots t_m^{n_m}.
\]

Tracing the proof of Proposition \ref{Taylor}, it is easy to see
that if the multihomogeneous ideal $I$ of $K\{X\}$ is invariant
under the formal partial derivatives, then the algebra of
constants $R_0$ is also multigraded. The following statement is an
analogue of the main theorem of \cite{D1}, see also \cite{D3},
Theorem 4.3.12. Our proof follows the idea of the proof given in
\cite{D3}.

\begin{corollary}\label{Hilbert}
Let $X=\{x_1,\ldots,x_m\}$ and let the ideal $I$ of $K\{X\}$ be
invariant under the formal derivatives. Then the Hilbert series of
the algebra $R=K\{X\}/I$ and its subalgebra of constants $R_0$ are
related by
\[
\text{\rm Hilb}(R,t_1,\ldots,t_m)= \prod_{j=1}^m\frac{1}{1-t_j}
\text{\rm Hilb}(R_0,t_1,\ldots,t_m).
\]
\end{corollary}

\begin{proof}
We fix a multihomogeneous basis $\{r_n\mid n=1,2,\ldots\}$ of
$R_0$. Let $\rho_j$ be the right multiplication by $x_j$,
$j=1,\ldots,m$. Then the Taylor expansion of Proposition
\ref{Taylor} shows that $R$ has a basis
\[
r_n\rho_1^{a_1}\cdots \rho_m^{a_m},\ n=1,2,\ldots,\ a_j\geq 0,\
j=1,\ldots,m.
\]
Hence, as a multigraded vector space $R$ is isomorphic to the
tensor product $R_0\otimes_KK[x_1,\ldots,x_m]$ with isomorphism
defined by
\[
r_n\rho_1^{a_1}\cdots \rho_m^{a_m}\to r_n\otimes x_1^{a_1}\cdots
x_m^{a_m}.
\]
Since the Hilbert series of the tensor product is equal to the
product of the Hilbert series of the factors, and the Hilbert
series of the polynomial algebra $K[x_1,\ldots,x_m]$ is
\[
\text{\rm Hilb}(K[x_1,\ldots,x_m],t_1,\ldots,t_m)
=\prod_{j=1}^m\frac{1}{1-t_j},
\]
we obtain immediately
\[
\text{\rm Hilb}(R,t_1,\ldots,t_m)= \text{\rm
Hilb}(R_0,t_1,\ldots,t_m) \text{\rm
Hilb}(K[x_1,\ldots,x_m],t_1,\ldots,t_m)
\]
\[
=\prod_{j=1}^m\frac{1}{1-t_j} \text{\rm Hilb}(R_0,t_1,\ldots,t_m).
\]
\end{proof}

\begin{remark}\label{youngmethod}

For $X=\{x_1,\ldots,x_m\}$, there is a natural action of the
general linear group $GL_m(K)$ on $K\{X\}$. More exactly, we can
consider the vector space $K\{X\}^{(k)}$ of the
elements of total degree $k$ as a finite dimensional representation of $GL_m(K)$.
Furthermore, for $n\leq m$, the symmetric group $S_n$ is embedded
into $GL_m(K)$ (permuting the first $n$ variables and fixing the
other variables). For the representations of symmetric and general
linear groups, see \cite{JK}, \cite{W}.

The irreducible representations of $S_n$ and $GL_m(K)$ are
symbolized by Young diagrams.

\medskip

 Under the conditions of Corollary \ref{Hilbert}, let us assume
that the ideal $I$ of $K\{X\}$ is invariant with respect to the
action of $GL_m(K)$. Thus $R^{(k)}$ is a representation of
$GL_m(K)$ for each $k=0,1,2,\ldots$~.

The proof of Corollary \ref{Hilbert} shows then that
\[
R^{(k)}=\bigoplus_{j=0}^k R_0^{(j)}\otimes \underbrace{\rm
 \vbox{\offinterlineskip \hbox{\boxtext{\phantom{X
}}\boxtext{\phantom{X }}\boxtext{$\cdots$}\boxtext{\phantom{X
}}\boxtext{\phantom{X }}}
 } \rm}_{k-j}
\]
as modules of the general linear group, where for simplicity
of notation we have used the Young diagram instead of the
corresponding $GL_m(K)$-module. An analogous
formula (but using induced modules in the place of tensor products)
holds if we look only at the multilinear parts ($n=m$) for
representations of the symmetric groups $S_n$.
\end{remark}

\begin{examples}
(i) Let $R=K\langle X\rangle$ be the free associative algebra. It
is a homomorphic image of $K\{X\}$ modulo the ideal generated by
the associators $(u,v,w)=(uv)w-u(vw)$, where $u,v,w$ run on the
set of all non-associative words in the alphabet $X$. Since
\[
\frac{\partial (u,v,w)}{\partial x_k}=
\left(\frac{\partial u}{\partial x_k},v,w\right)+
\left(u,\frac{\partial v}{\partial x_k},w\right)+
\left(u,v,\frac{\partial w}{\partial x_k}\right),
\]
this ideal is invariant under all partial derivatives.
The algebra of constants of $K\langle X\rangle$
is spanned by all products of commutators
$[x_{j_1},\ldots ,x_{j_p}]\cdots[x_{k_1},\ldots,x_{k_q}]$,
where $p,\ldots,q\geq 2$. See the comments in the introduction of the
present
paper as well as \cite{D3}, Section 4.3.
\par
(ii) Let $\mathfrak M$ be a variety of unitary (associative or
non-associative) algebras. The T-ideal $T({\mathfrak M})\subset
K\{x_1,x_2,\ldots\}$ of all polynomial identities of $\mathfrak M$
is generated by constants. For associative algebras this is well
known, see e.~g. \cite{D3}, Section 4.3 for the proof. For
non-associative algebras see \cite{D2}, Corollary 1.6.
\par
(iii) Let $L=L(X)$ be the free Lie algebra freely generated by the
set $X$. By the Witt theorem we may assume that $L(X)$ is the Lie
subalgebra of $K\langle X\rangle$ generated by $X$ with respect to
the bracket multiplication $[u,v]=uv-vu$. Let $I_L$ be a
multigraded ideal of $L(X)$ and let $G=L(X)/I_L$ be the
corresponding multigraded Lie algebra. We assume that $I_L$ does
not contain elements of first degree, i.~e. $I_L$ is contained in
the commutator ideal $L'(X)$ of $L(X)$. Clearly, the higher
commutators $u=[x_{k_1},\ldots,x_{k_n}]$ which span $L'(X)$ vanish
under formal partial derivatives. By the Poincar\'e-Birkhoff-Witt
theorem, the universal enveloping algebra of $G$ is the
associative algebra $U(G)=K\langle X\rangle/I$, where $I$ is the
ideal of $K\langle X\rangle$ generated by $I_L$. Hence we may
apply Proposition \ref{Taylor}. A concrete basis of the constants
in $U(G)$ can be obtained in the following way. Let us fix a basis
of $G$ consisting of $X$ and some higher commutators
$u_j=[x_{k_{j1}},\ldots,x_{k_{jp_j}}]$, $p_j\geq 2$,
$j=1,2,\ldots$. Then $U(G)$ has a basis
\[
\{u_1^{b_1}\cdots u_n^{b_n}x_1^{a_1}\cdots x_m^{a_m}\mid a_j,b_k\geq
0\}
\]
and the algebra $U_0=(U(G))_0$ of the constants of $U(G)$ is spanned by
the
basis elements with $a_1=\cdots=a_m=0$.
\par
(iv) A special case of (iii) is the free metabelian Lie algebra
\begin{equation*}
M(X)=L(X)/L''(X).
\end{equation*}
It has a basis
\[
X\cup\{[x_{k_1},x_{k_2},\ldots,x_{k_n}]\mid k_1>k_2\leq \cdots\leq
k_n\}
\]
and the higher commutators commute.

 Hence the algebra of constants
of $U(M(X))$ is isomorphic to the polynomial algebra generated by
the ``commutator variables'' $[x_{k_1},x_{k_2},\ldots,x_{k_n}]$,
where $k_1>k_2\leq \cdots\leq k_n$. In \cite{G3}, Gerritzen
used the completion of the algebra $U(M(x,y))$ in order to find a
simple expression for the evaluation, modulo the metabelian
identity $[[x_1,x_2],[x_3,x_4]]=0$, of the Hausdorff series
$z=x+y+[x,y]/2+\cdots$, where
$\text{exp}(z)=\text{exp}(x)\text{exp}(y)$.
\end{examples}

\begin{remark} (cf.\ \cite{GH} Proposition 2.3.)
Let $X=\{x_1,\ldots,x_m\}$ be any finite set and let the
multihomogeneous ideal $I$ of $K\{X\}$ be invariant under the
formal partial derivatives. Let $R=K\{X\}/I$ be the factor
algebra. The augmentation ideal $\omega(R)$ of $R$ is generated by
the set $X$ and consists of all elements of $R$ without constant
terms. The powers of $\omega(R)$ induce a topology on $R$ called
the formal power series topology. It is defined by the property
that the sets $r+\omega^n(R)$, $r\in R$, $n\geq 1$, form a basis
for the open sets. Clearly, the completion $\hat R$ of $R$ with
respect to this topology consists of all formal power series in
$X$, i.~e. formal infinite sums of the form
\[
r=\sum_{n\geq 0}r_n,
\]
where $r_n\in R$ are homogeneous elements of degree $n$.
Similarly, the completion $\widehat{R_0}$ of the subalgebra $R_0$
of the constants in $R$ consists of the formal power series
$c=\sum_{k\geq 0}c_k$, where $c_k$ are homogeneous constants of
degree $k$, $k\geq 0$. The formal partial derivatives can be
extended in a unique continuous way to derivatives of $\hat R$ and
it is clear that the algebra of constants of $\hat R$ coincides
with $\widehat{R_0}$.

Using the Taylor formula, every homogeneous $r_n\in R$ can be
expressed in the form
\[
r_n=\sum_{k\geq 0}c_{nk}\rho_1^{p_1}\cdots\rho_m^{p_m},
\]
where $c_{nk}\in R_0$ are homogeneous constants of degree $k$,
$\rho_j$ is the operator of the right multiplication by $x_j$ and
$p_j\geq 0$, $j=1,\ldots,m$, $p_1+\cdots+p_m=n-k$. Hence we get a
Taylor expansion
\[
r=\sum_{p_j\geq 0}\left(\sum_{k\geq
0}c_{nk}\right)\rho_1^{p_1}\cdots\rho_m^{p_m}
=\sum_{p_j\geq 0}c_p\rho_1^{p_1}\cdots\rho_m^{p_m},\
c_p\in\widehat{R_0},
\]
for every $r\in \hat R$.

\end{remark}

\begin{theorem}\label{general}

Let $R=K\{X\}/I$, where the multihomogeneous ideal $I$ is
invariant under the formal derivatives. For each $x_j\in X$ and
each $k\geq 1$ we fix an operator $\mu_{jk}=\mu_{jk}(x_j)$ in the
algebra of multiplications such that $\mu_{jk}$ is homogeneous of
degree $k$, depends on $x_j$ only and satisfies the property
$\mu_{jk}(1)\not=0$.

Then every element $r$ of $R$ has a presentation of the form
\[
r=\sum r_a\mu_{1a_1}\cdots\mu_{ma_m},\
r_a\in R_0,\ a_1,\ldots,a_m\geq 0,
\]
and this presentation is unique.
\end{theorem}

\begin{proof}
Since the algebra $R$ is multihomogeneous, it is sufficient to
consider the subalgebra generated by $x_1,\ldots,x_m$, i.~e. we may
assume that
$X=\{x_1,\ldots,x_m\}$.

We have that
$\mu_{jk}(x_j)=\sum\alpha_{\nu,u}\nu_{1,u_1}\cdots\nu_{n,u_n}$,
where $\alpha_{\nu,u}\in K$ and $\nu_{su_s}=\lambda_{u_s}$ or
$\nu_{su_s}=\rho_{u_s}$ for some monomial $u_s=u_s(x_j)$. The
condition $\mu_{jk}(1)\not=0$ just means
$\sum\alpha_{\nu,u}\not=0$.

It is easy to see that for $r_0\in R_0$
\[
\frac{\partial^k r_0\mu_{1k}(x_1)}{\partial x_1^k}=k!r_0\mu_{1k}(1).
\]
We fix a multihomogeneous basis $\{r_n\mid n=1,2,\ldots\}$ of $R_0$.
\par
First, we want to show that the elements
$r_n\mu_{1a_1}\cdots\mu_{ma_m}$, $n\geq 1$, $a_j\geq 0$, are
linearly independent. Let
\[
r=\sum \beta_{n,a}r_n\mu_{1a_1}\cdots\mu_{ma_m}=0,
\]
where some of the $\beta_{n,a}\in K$ are nonzero.

 Choosing the
maximal $m$-tuple $(k_1,\ldots,k_m)$ among all indices
$(a_1,\ldots,a_m)$ with $\beta_{n,a}\not=0$, with respect to the
lexicographical order, we obtain for $c=k_1+\cdots+k_m$
\[
\frac{\partial^cr}{\partial x_1^{k_1}\cdots\partial x_m^{k_m}}
=k_1!\cdots k_m!\mu_{1k_1}(1)\cdots\mu_{mk_m}(1)\sum \beta_{n,k}r_n=0.
\]
Hence $\beta_{n,k}=0$ for all $n$ which is a contradiction.
\par
Now, let $S$ be the vector subspace of $R$ spanned by all
$r_n\mu_{1a_1}\cdots\mu_{ma_m}$, $n\geq 1$, $a_j\geq 0$. Since these
elements
are linearly independent, they form a multigraded basis of $S$. Bearing
in mind that $\text{\rm deg }\mu_{jk}=k$,
we obtain for the Hilbert series of $S$ and $R$
\[
\text{\rm Hilb}(S,t_1,\ldots,t_m)=\text{\rm Hilb}(R_0,t_1,\ldots,t_m)
\sum t_1^{a_1}\cdots t_m^{a_m}
\]
\[
=\text{\rm Hilb}(R_0,t_1,\ldots,t_m)
\text{\rm Hilb}(K[x_1,\ldots,x_m],t_1,\ldots,t_m)
=\text{\rm Hilb}(R,t_1,\ldots,t_m).
\]
The coincidence of the Hilbert series of $S$ and $R$ gives that $S=R$.
\end{proof}

\begin{example}
Let $SJ(X)$ be the free special Jordan algebra generated by the set
$X$. This is the Jordan subalgebra of $K\langle X\rangle$
generated by $X$ with respect to the multiplication $u\circ
v=uv+vu$. In $K\langle X\rangle$ we define the operators
$\mu_{jk}=(\lambda_{x_j}+\rho_{x_j})^k$. Clearly
$u\mu_{jk}=(\cdots(u\circ x_j)\circ\cdots)\circ x_j$ ($k$-times
Jordan multiplication by $x_j$). If $SJ_0(X)$ is the subalgebra of
constants, then Theorem \ref{general} gives that the elements of
$SJ(X)$ have the form $\sum (\cdots(r_j\circ
x_{j_1})\circ\cdots)\circ x_{j_n}$, where $x_{j_1}\leq\cdots\leq
x_{j_n}$ and this presentation is unique.
\end{example}

\section{Constants in free non-associative algebras}

A variety $\mathfrak M$ of algebras satisfies the Nielsen-Schreier
property if any subalgebra of the relatively free algebra
$F({\mathfrak M})=K\{X\}/T({\mathfrak M})$ is again relatively
free in $\mathfrak M$. Examples of such varieties are the class of
all algebras \cite{K}, all commutative (non-associative) algebras
and all anticommutative algebras, Shirshov \cite{Sh}, all Lie
algebras, etc. See the paper by Umirbaev \cite{U} for the
description of the varieties with Nielsen-Schreier property in
terms of necessary and sufficient conditions. (Of course, the free
anticommutative algebra and the free Lie algebra do not satisfy
our condition for invariantness under the formal partial
derivatives because they are not unitary.) If the condition for
invariantness under the formal partial derivatives is satisfied,
we may consider the problem for the description of the algebra of
constants and its free generating set.

\medskip

Let $K\{x\}$ be the free non-associative algebra generated by a
single element $x$. Recall that the number of all non-associative
and noncommutative words in $x$ of length $n$ is equal to the
$n$-th Catalan number
\[
c_n=\frac{1}{n}\binom{2n-2}{n-1}
\]
The sequence $c_n$, $n=1,2,\ldots$, satisfies the relation
\[
c_n=\sum_{p=1}^{n-1}c_pc_{n-p},\ n=2,3,\ldots
\]
and the first few Catalan numbers are
\[
c_1=c_2=1,c_3=2,c_4=5,\ldots \ .
\]
If we define $c_0=1$, then the generating function of the sequence
is
\[
c(t)=\sum_{n\geq 0}c_nt^n=1+\frac{1-\sqrt{1-4t}}{2}.
\]
(It is also convenient to define $c_{-1}=0$.)

\bigskip

\begin{lemma}\label{generating}
{\rm (See e.~g. \cite{G2}, Proposition 1.1.)} Let $X$ be a
(multi)graded set with generating function $g(T)\in
\mathbb{Z}[[T]]$, where $T=\{t_1,t_2,\ldots\}$ is some set of
commuting variables. Then the Hilbert series of the absolutely
free algebra $K\{X\}$ is
\[
\text{\rm Hilb}(K\{X\},T)=c(g(T))=1+\frac{1-\sqrt{1-4g(T)}}{2}.
\]
\end{lemma}

\begin{proof}
The considerations are standard and similar to those in \cite{G2}.
We include them for completeness of the exposition only. Let
$\{x\}$ be the set of all non-associative noncommutative words in
one variable. If we replace any word $u=(x\cdots x)(x\cdots x)\in
\{x\}$ of length $n$ with $(x_{j_1}\cdots
x_{j_p})(x_{j_{p+1}}\cdots x_{j_n})$, where the $x_j$'s run on
$X$, we obtain a basis of $K\{X\}$. Hence the generating function
of the words of length $n$ in $X$ is equal to $c_ng^n(T)$ and the
Hilbert series of $K\{X\}$ is
\[
\text{\rm Hilb}(K\{X\},T)=\sum_{n\geq 0}c_ng^n(T)=c(g(T))
=1+\frac{1-\sqrt{1-4g(T)}}{2}.
\]
\end{proof}

\begin{remark}

Let $C=(K\{x\})_0$ be the algebra of constants of $K\{x\}$.

To determine $\gamma(t)=\text{\rm Hilb}(C,t)=\sum_{n\geq
0}\gamma_nt^n$, we can use Corollary \ref{Hilbert} to get that
\[
\text{\rm Hilb}(K\{x\},t)=\frac{1}{1-t}\text{\rm Hilb}(C,t)=
\left(1+\frac{1-\sqrt{1-4t}}{2}\right).
\]

A theorem of Kurosh \cite{K} states that for any generating set
$X$ any subalgebra $A$ of the algebra $K\{X\}$ is free. Of course,
if $A$ is graded with respect to any grading of $K\{X\}$, then the
set of free generators of $A$ is also graded. Let $Y$ be a
homogeneous system of free generators of $C=K\{Y\}$.

Lemma \ref{generating} and direct calculations show that
\[
\gamma(t)=(1-t)\left(1+\frac{1-\sqrt{1-4t}}{2}\right)
=1+\frac{1-\sqrt{1-4g(t)}}{2},
\]
\[
g(t)=-3t^2+t^3+3t(1-t)\frac{1-\sqrt{1-4t}}{2} =t^3+3\sum_{n\geq
4}(c_{n-1}-c_{n-2})t^n.
\]

Thus we have obtained a different proof of \cite{GH} Proposition
2.9.\ and \cite{GH} Proposition 2.10.(ii):

\end{remark}

\begin{proposition}\label{constants}(cf.\ \cite{GH}, Prop.\ 2.9.\ and 2.10.)
Let $C=(K\{x\})_0$ be the algebra of constants of $K\{x\}$ with
respect to the formal derivative $d/dx$.
\par
{\rm (i)} The Hilbert series of $C$ is
\[
\text{\rm Hilb}(C,t)=(1-t)\left(1+\frac{1-\sqrt{1-4t}}{2}\right)
=\sum_{n\geq 0}(c_n-c_{n-1})t^n.
\]
\par
{\rm (ii)} Let $Y$ be a homogeneous system of free generators of
$C=K\{Y\}$. Then the generating function of $Y$ is
\[
g(t)=-3t^2+t^3+3t(1-t)\frac{1-\sqrt{1-4t}}{2} =t^3+3\sum_{n\geq
4}(c_{n-1}-c_{n-2})t^n.
\]
\end{proposition}
\qed

\begin{remark}

A linear ordering $\prec$ on the set $\{x\}$ of words can be
defined by:

\begin{enumerate}

\item[(i)] If $u,v\in \{x\}$ and $\text{\rm deg }u<\text{\rm deg }v$, then
$u\prec v$.

\item[(ii)] If $\text{\rm deg }u=\text{\rm deg }v>1$, where $u=u_1u_2$,
$v=v_1v_2$ and $u_2\prec v_2$, then $u\prec v$.

\item[(iii)] If $\text{\rm deg }u=\text{\rm deg }v>1$, where $u=u_1w$,
$v=v_1w$ and $u_1\prec v_1$, then $u\prec v$.

\end{enumerate}

For example,
\[
x\prec x^2\prec (xx)x\prec x(xx)\prec
\]
\[
\prec ((xx)x)x \prec (x(xx))x \prec (xx)(xx)\prec x((xx)x)\prec
x(x(xx))\prec\cdots \ .
\]
For any nonzero element
\[
f=\sum_{k=1}^p\alpha_ku_k\in K\{x\},\ 0\not=\alpha_k\in K,\ u_k\in
\{x\},\ u_1\prec \cdots\prec u_p,
\]
we define the leading term of $f$ as $\text{\rm
lt}(f)=\alpha_pu_p$. Clearly, if $f$ and $g$ are two nonzero
polynomials in $K\{x\}$, then $\text{\rm lt}(fg)=\text{\rm
lt}(f)\text{\rm lt}(g)$.

Using this ordering one can show, see \cite{GH} Proposition 2.9.,
that $C$ has a vector space basis consisting of all polynomials
\[
u-\frac{u'}{1!}+\frac{u''}{2!}-\frac{u'''}{3!}+\cdots,
\]
where $u$ runs on the set of all words in $\{x\}$ which are not in
the form $u=v\cdot x$, $v\in\{x\}$.

\medskip

 One of the basic properties of the free
algebras in any variety $\mathfrak M$ satisfying the
Nielsen-Schreier property is the following. If we define an
ordering which agrees with the multiplication in $F({\mathfrak
M})$, and if some nonzero polynomials $f_1,\ldots,f_m\in
F({\mathfrak M})$ satisfy a nontrivial relation
$\omega(f_1,\ldots,f_m)=0$, then the leading term $\text{\rm
lt}(f_i)$ of one of the polynomials belongs to the subalgebra
generated by the other leading terms $\text{\rm
lt}(f_1),\ldots,\text{\rm lt}(f_{i-1}), \text{\rm
lt}(f_{i+1}),\ldots,\text{\rm lt}(f_m)$.

Applying this argument in the case of $K\{x\}$ one can deduce that
a free algebra basis of $C$ is given by the set $Y$ consisting of
those
\[
u-\frac{u'}{1!}+\frac{u''}{2!}-\frac{u'''}{3!}+\cdots,
\]
which have the property that the word $u$ has one of the forms
\[
u=v\cdot (x^2),\ u=x\cdot w,\ u=(x^2)\cdot w,\ u=(v_1\cdot x)w,\
u=v(w_1\cdot x),
\]
where $\text{\rm deg }v,\text{\rm deg }w,\text{\rm deg
}v_1,\text{\rm deg }w_1\geq 2$; see \cite{GH} Proposition
2.10.(i).

\end{remark}

\begin{remark}\label{constmag}

Let now $X=\{x_1,\ldots,x_m\}$ and $C=(K\{X\})_0$ be the algebra
of constants of $K\{X\}$. We have already observed, see
 Remark \ref{youngmethod}, that
\[
K\{X\}^{(k)}=\bigoplus_{j=0}^k C^{(j)}\otimes \underbrace{\rm
 \vbox{\offinterlineskip \hbox{\boxtext{\phantom{X
}}\boxtext{\phantom{X }}\boxtext{$\cdots$}\boxtext{\phantom{X
}}\boxtext{\phantom{X }}}
 } \rm}_{k-j}
\]
as modules of the general linear group. (Note that the
conditions given in Remark \ref{youngmethod} are trivially
fulfilled.)

The representation given by the component $K\langle
X\rangle^{(k)}$ of degree $k$ of the free associative algebra is well-known, and
$K\{X\}^{(k)}$ consists just of $c_k$ copies.

We can apply methods described by Regev (cf.\ \cite{R}), and of
\cite{D1}, \cite{D3}, to recursively determine $C^{(k)}$ for
$k=0,1,2,\ldots$, starting with $C^{(0)}=K$, $C^{(1)}=0$. The
methods make use of the Young rule and Littlewood-Richardson rule.

\end{remark}

\begin{example} From
$K\{X\}^{(2)}=K\langle X\rangle^{(2)}= {\rm
 \vbox{\offinterlineskip \hbox{\boxtext{\phantom{X
}}\boxtext{\phantom{X }}}
 } \rm}\oplus C^{(2)}$ we get
 that $C^{(2)}$ corresponds to the sign representation ${\rm
 \vbox{\offinterlineskip \hbox{\boxtext{\phantom{X
}}} \hbox{\boxtext{\phantom{X }}}
 } \rm}\ ,$ provided $m\geq 2$. (The generators are the commutators
 $[x_i,x_j]$.)
 In the following, let $m\geq k$ always.

In degree 3 we have to compare $c_2=2$ copies of $K\langle
X\rangle^{(3)}$ with

\[
 {\rm \vbox{\offinterlineskip
\hbox{\boxtext{\phantom{X }}\boxtext{\phantom{X
}}\boxtext{\phantom{X }}}
 } \rm}\oplus\Bigl(\
{\rm
 \vbox{\offinterlineskip \hbox{\boxtext{\phantom{X
}}} \hbox{\boxtext{\phantom{X }}}
 } \rm}\otimes \boxtext{\phantom{X }}\ \Bigr)\oplus\ C^{(3)}.
\]
By Young's rule,

\[
{\rm
 \vbox{\offinterlineskip \hbox{\boxtext{\phantom{X
}}} \hbox{\boxtext{\phantom{X }}}
 } \rm}\otimes \boxtext{\phantom{X }}\
=\ {\rm
 \vbox{\offinterlineskip \hbox{\boxtext{\phantom{X
}}\boxtext{\phantom{X }}} \hbox{\boxtext{\phantom{X }}}
 } \rm}
\oplus\  {\rm
 \vbox{\offinterlineskip \hbox{\boxtext{\phantom{X
}}} \hbox{\boxtext{\phantom{X }}} \hbox{\boxtext{\phantom{X }}} }
\rm}\ .
\]
Therefore

\[
C^{(3)} =
 {\rm \vbox{\offinterlineskip
\hbox{\boxtext{\phantom{X }}\boxtext{\phantom{X
}}\boxtext{\phantom{X }}}
 } \rm}\oplus\ 3
\ {\rm
 \vbox{\offinterlineskip \hbox{\boxtext{\phantom{X
}}\boxtext{\phantom{X }}} \hbox{\boxtext{\phantom{X }}}
 } \rm}
\oplus\  {\rm
 \vbox{\offinterlineskip \hbox{\boxtext{\phantom{X
}}} \hbox{\boxtext{\phantom{X }}} \hbox{\boxtext{\phantom{X }}} }
\rm}\ .
\]

\medskip

It is easily seen that the multiplicity of the trivial
representation in $C^{(k)}$ is always $c_k-c_{k-1}$ (we already
know this from the case $X=\{x\}$).

For $k=4$, $c_k$ is 5, and
\[
\bigoplus_{j=0}^3 C^{(j)}\otimes \underbrace{\rm
 \vbox{\offinterlineskip \hbox{\boxtext{\phantom{X
}}\boxtext{$\cdots$}\boxtext{\phantom{X }}}
 } \rm}_{4-j}
\]
is computed using that
\[
{\rm
 \vbox{\offinterlineskip \hbox{\boxtext{\phantom{X
}}} \hbox{\boxtext{\phantom{X }}}
 } \rm}\otimes\
{\rm
 \vbox{\offinterlineskip \hbox{\boxtext{\phantom{X
}}\boxtext{\phantom{X }}}
 } \rm}
\ = \ {\rm
 \vbox{\offinterlineskip \hbox{\boxtext{\phantom{X }}\boxtext{\phantom{X
}}\boxtext{\phantom{X }}} \hbox{\boxtext{\phantom{X }}}
 } \rm}
 \oplus\ {\rm
 \vbox{\offinterlineskip \hbox{\boxtext{\phantom{X
}}\boxtext{\phantom{X }}} \hbox{\boxtext{\phantom{X }}}
\hbox{\boxtext{\phantom{X }}} } \rm}\ .
\]

Then one gets that $C^{(4)}$ is given by
\[
 3 {\rm \vbox{\offinterlineskip
\hbox{ \boxtext{\phantom{X }}\boxtext{\phantom{X
}}\boxtext{\phantom{X }}\boxtext{\phantom{X }}}
 } \rm}\ \oplus\ 10\
{\rm
 \vbox{\offinterlineskip \hbox{\boxtext{\phantom{X
}}\boxtext{\phantom{X }}\boxtext{\phantom{X }}}
\hbox{\boxtext{\phantom{X }}}
 } \rm}
\oplus\ 7\ {\rm
 \vbox{\offinterlineskip \hbox{\boxtext{\phantom{X
}}\boxtext{\phantom{X }}} \hbox{\boxtext{\phantom{X
}}\boxtext{\phantom{X }}}
 } \rm}
 \oplus\ 10\ {\rm
 \vbox{\offinterlineskip \hbox{\boxtext{\phantom{X
}}\boxtext{\phantom{X }}} \hbox{\boxtext{\phantom{X }}}
\hbox{\boxtext{\phantom{X }}} } \rm}
 \oplus\ 4\ {\rm
 \vbox{\offinterlineskip \hbox{\boxtext{\phantom{X
}}} \hbox{\boxtext{\phantom{X }}} \hbox{\boxtext{\phantom{X }}}
\hbox{\boxtext{\phantom{X }}} } \rm} \ .
\]

\end{example}

\bigskip

\begin{remark}
The free commutative non-associative algebra $R$ in $m$ variables
has a nice basis and allows the formal derivatives $d/dx_i$. Its
Hilbert series $\mathrm{Hilb}(R,t)$ given by
\[
\mathrm{Hilb}(R,t)=mt+\frac{1}{2}\mathrm{Hilb}(R,t)^2
+\frac{1}{2}\mathrm{Hilb}(R,t^2)
\]
can be written in the form
\begin{equation*}
\begin{split}
\mathrm{Hilb}(R,t)&=1-\sqrt{1-\mathrm{Hilb}(R,t^2)-2mt}\\
&=1-\sqrt{\sqrt{1-\mathrm{Hilb}(R,t^4)-2mt^2}-2mt}=\cdots\\
\end{split}
\end{equation*}
(see \cite{P}).

Using the techniques described above it is possible to determine
the first components of $R_0$.

\end{remark}

\begin{example}
We can recursively compute $R_0^{(k)}$ for the free commutative
non-associative algebra $R$ generated by $m$ variables, assuming
$m\geq k$ always. Clearly $R_0^{(0)}=K$ and $R_0^{(1)}=0$.

Since
\[
R_0^{(0)}\otimes {\rm
 \vbox{\offinterlineskip \hbox{\boxtext{\phantom{X
}}\boxtext{\phantom{X }}}
 } \rm}
\]
is already $R^{(2)}$, we get that  $R_0^{(2)}=0$, i.~e.\ there are
no constants in degree 2.

In degree 3, we observe that $R^{(3)}$ is
\[
 {\rm \vbox{\offinterlineskip
\hbox{\boxtext{\phantom{X }}\boxtext{\phantom{X
}}\boxtext{\phantom{X }}}
 } \rm}\oplus
{\rm
 \vbox{\offinterlineskip \hbox{\boxtext{\phantom{X
}}\boxtext{\phantom{X }}} \hbox{\boxtext{\phantom{X }}}
 } \rm}
\]
and the second summand must be equal to $R_0^{(3)}$.

Now
\[
\bigoplus_{j=0}^3 R_0^{(j)}\otimes \underbrace{\rm
 \vbox{\offinterlineskip \hbox{\boxtext{\phantom{X
}}\boxtext{$\cdots$}\boxtext{\phantom{X }}}
 } \rm}_{4-j}
\]
\[
 = {\rm \vbox{\offinterlineskip
\hbox{\boxtext{\phantom{X }}\boxtext{\phantom{X
}}\boxtext{\phantom{X }}\boxtext{\phantom{X }}}
 } \rm}\oplus
{\rm
 \vbox{\offinterlineskip \hbox{\boxtext{\phantom{X
}}\boxtext{\phantom{X }}\boxtext{\phantom{X }}}
\hbox{\boxtext{\phantom{X }}}
 } \rm}
\oplus\ {\rm
 \vbox{\offinterlineskip \hbox{\boxtext{\phantom{X
}}\boxtext{\phantom{X }}} \hbox{\boxtext{\phantom{X
}}\boxtext{\phantom{X }}}
 } \rm}
 \oplus\ {\rm
 \vbox{\offinterlineskip \hbox{\boxtext{\phantom{X
}}\boxtext{\phantom{X }}} \hbox{\boxtext{\phantom{X }}}
\hbox{\boxtext{\phantom{X }}} } \rm}
\]

as
\[
{\rm
 \vbox{\offinterlineskip \hbox{\boxtext{\phantom{X
}}\boxtext{\phantom{X }}} \hbox{\boxtext{\phantom{X }}}
 } \rm}\otimes\ \boxtext{\phantom{X }} =
{\rm
 \vbox{\offinterlineskip \hbox{\boxtext{\phantom{X
}}\boxtext{\phantom{X }}\boxtext{\phantom{X }}}
\hbox{\boxtext{\phantom{X }}}
 } \rm}
\oplus\ {\rm
 \vbox{\offinterlineskip \hbox{\boxtext{\phantom{X
}}\boxtext{\phantom{X }}} \hbox{\boxtext{\phantom{X
}}\boxtext{\phantom{X }}}
 } \rm}
 \oplus\ {\rm
 \vbox{\offinterlineskip \hbox{\boxtext{\phantom{X
}}\boxtext{\phantom{X }}} \hbox{\boxtext{\phantom{X }}}
\hbox{\boxtext{\phantom{X }}} } \rm}
\]
by Young's rule.

Since $R^{(4)}$ is given by

\[
  2 {\rm \vbox{\offinterlineskip
\hbox{ \boxtext{\phantom{X }}\boxtext{\phantom{X
}}\boxtext{\phantom{X }}\boxtext{\phantom{X }}}
 } \rm}\ \oplus\ 2\
{\rm
 \vbox{\offinterlineskip \hbox{\boxtext{\phantom{X
}}\boxtext{\phantom{X }}\boxtext{\phantom{X }}}
\hbox{\boxtext{\phantom{X }}}
 } \rm}
\oplus\ 2\ {\rm
 \vbox{\offinterlineskip \hbox{\boxtext{\phantom{X
}}\boxtext{\phantom{X }}} \hbox{\boxtext{\phantom{X
}}\boxtext{\phantom{X }}}
 } \rm}
 \oplus\ 1\ {\rm
 \vbox{\offinterlineskip \hbox{\boxtext{\phantom{X
}}\boxtext{\phantom{X }}} \hbox{\boxtext{\phantom{X }}}
\hbox{\boxtext{\phantom{X }}} } \rm}
\]

we get that
 $R_0^{(4)}$ must be equal to
\[
  {\rm \vbox{\offinterlineskip
\hbox{\boxtext{\phantom{X }}\boxtext{\phantom{X
}}\boxtext{\phantom{X }}\boxtext{\phantom{X }}}
 } \rm}\oplus
{\rm
 \vbox{\offinterlineskip \hbox{\boxtext{\phantom{X
}}\boxtext{\phantom{X }}\boxtext{\phantom{X }}}
\hbox{\boxtext{\phantom{X }}}
 } \rm}
\oplus\ {\rm
 \vbox{\offinterlineskip \hbox{\boxtext{\phantom{X
}}\boxtext{\phantom{X }}} \hbox{\boxtext{\phantom{X
}}\boxtext{\phantom{X }}}
 } \rm}\ .
\]

\end{example}

\section{Linear Differential Equations}

We consider ordinary linear differential equations of the form
 $$
y^{(n)}+a_1y^{(n-1)}+\cdots+a_{n-1}y'+a_ny=f(x), $$
 where the
coefficients $a_1,\ldots,a_{n-1},a_n$ are constants from the base
field $K$ and $f(x)$ is a formal power series from the completion
$\hat R$ of the algebra $R=K\{x\}/I$ where $I$ is a homogeneous
ideal invariant under the formal derivative $d_x=d/dx$.

\begin{proposition}\label{linear equations}
Let the homogeneous ideal $I$ of $K\{x\}$ be invariant under the formal
derivative
$d/dx$ and $R=K\{x\}/I$.
For any constants $c_0(x),c_1(x),\ldots,c_{n-1}(x)$
in the algebra $\hat R$ of formal power series the linear differential
equation
$$
y^{(n)}+a_1y^{(n-1)}+\cdots+a_{n-1}y'+a_ny=f(x),
$$
$f(x)\in \hat R$, $a_1,\ldots,a_{n-1},a_n\in K$,
has a unique solution $y(x)$ of the form
\[
y(x)=c_0+c_1\frac{\rho}{1!}+c_2\frac{\rho^2}{2!}+\cdots
+c_{n-1}\frac{\rho^{n-1}}{(n-1)!}+c_n\frac{\rho^n}{n!}+\cdots,
\]
where $\rho$ is the operator of right multiplication by $x$ and
$c_n(x),c_{n+1}(x),\ldots$ belong to $\widehat{R_0}$.
\end{proposition}

\begin{proof}
Let $f(x)=f_0+f_1\rho/1!+f_2\rho^2/2!+\cdots$, where
$f_k\in\widehat{R_0}$.
We are looking for a solution of the form
$y(x)=c_0+c_1\rho/1!+c_2\rho^2/2!+\cdots$, $c_k\in\widehat{R_0}$,
where the first $n$ coefficients $c_0,c_1,\ldots,c_{n-1}$ coincide with
the
prescribed ones. Since $c_k'=0$, it is directly to see that
\[
y^{(j)}(x)=c_j+c_{j+1}\frac{\rho}{1!}+c_{j+2}\frac{\rho^2}{2!}+\cdots
\]
and the differential equation has the form
\[
\sum_{j\geq 0}(c_{j+n}+a_1c_{j+n-1}+\cdots+a_{n-1}c_{j+1}+a_nc_j)\rho^j
=\sum_{j\geq 0}f_j\rho^j.
\]
Comparing the coefficients of the power series in $\rho$, we obtain
that
\[
c_{j+n}+a_1c_{j+n-1}+\cdots+a_{n-1}c_{j+1}+a_nc_j=f_j.
\]
Since $c_0,c_1,\ldots,c_{n-1}$ are already fixed, and hence known, this
allows
to define step-by-step and in a unique way
the other coefficients $c_n,c_{n+1},\ldots$~.
\end{proof}

As in the case of functions in one real or complex variable,
all solutions of any ordinary linear differential equation
are obtained as sums of a given partial solution of the given equation
and all the solutions of the corresponding homogeneous equation.
In this setup, we establish an analogue of the well know result
about the general form of the solutions of the homogeneous equation.
For simplicity of the exposition, we assume that
$K=\mathbb{C}$.

\begin{theorem}\label{homogeneous case}
Let $a_1,\ldots,a_{n-1},a_n\in\mathbb{C}$ and let
$\lambda_1,\ldots,\lambda_p\in \mathbb{C}$
be all pairwise different solutions of the algebraic equation
\[
\lambda^n+a_1\lambda^{n-1}+\cdots+a_{n-1}\lambda+a_n=0
\]
with multiplicity $k_1,\ldots,k_p$, respectively.
Let the homogeneous ideal $I$ of $\mathbb{C}\{x\}$
be invariant under the formal derivative
$d/dx$ and let $R=\mathbb{C}\{x\}/I$.
Then all solutions of the homogeneous linear differential equation
$$
y^{(n)}+a_1y^{(n-1)}+\cdots+a_{n-1}y'+a_ny=0
$$
in the algebra $\hat R$ of formal power series
are given by the formula
\[
y(x)=(c_{10}+c_{11}\rho+\cdots+c_{1,k_1-1}\rho^{k_1-1})
\text{\rm exp}(\lambda_1\rho)+
\cdots
\]
\[
+(c_{p0}+c_{p1}\rho+\cdots+c_{p,k_p-1}\rho^{k_p-1})
\text{\rm exp}(\lambda_p\rho),
\]
where $c_{ij}$ are arbitrary constants in $\widehat{R_0}$,
$\rho$ is the operator of right multiplication by $x$ and
$\text{\rm exp}(\lambda\rho)$,
$\lambda\in\mathbb{C}$, is the element of the completion
$\widehat{\mathcal{M}(R)}$
of the algebra of multiplications $\mathcal{M}(R)$ defined by
\[
\text{\rm exp}(\lambda\rho)=1+\lambda\frac{\rho}{1!}
+\lambda^2\frac{\rho^2}{2!}+\lambda^3\frac{\rho^3}{3!}+\cdots
\]
\end{theorem}

\begin{proof}
If $c(x)\in\widehat{R_0}$ is any constant and
$s(\rho)=\sum_{j\geq 0}\alpha_j\rho^j$, $\alpha_j\in \mathbb{C}$,
is a formal power series in the operator of multiplication $\rho$,
then,
since $c'(x)=0$,
\[
\frac{dc(x)s(\rho)}{dx}=\sum_{j\geq 0}\alpha_j\frac{dc(x)\rho^j}{dx}=
c(x)\sum_{j\geq 0}\alpha_jj\rho^j=c(x)s'(\rho).
\]
Hence $y(x)=c(x)s(\rho)$ is a solution of the equation
$y^{(n)}+a_1y^{(n-1)}+\cdots+a_{n-1}y'+a_ny=0$ if and only if
$s(t)$ is a solution of the equation
$s^{(n)}+a_1s^{(n-1)}+\cdots+a_{n-1}s'+a_ns=0$.
The standard theory of homogeneous linear differential equations
with constant coefficients gives that the $n$ formal power series
\[
s_{ij}(t)=t^j\text{\rm exp}(\lambda_it),
\ j=0,1,\ldots,k_i-1,\
i=1,\ldots,p,
\]
are solutions of the  equation
$s^{(n)}(t)+a_1s^{(n-1)}(t)+\cdots+a_{n-1}s'(t)+a_ns(t)=0$
in $\mathbb{C}[[t]]$. Hence
\[
y_{ij}=y_{ij}(c_{ij},x)=c_{ij}\rho^j\text{\rm exp}(\lambda_i\rho),
\ j=0,1,\ldots,k_i-1, \
i=1,\ldots,p,
\]
are solutions of the equation
$y^{(n)}+a_1y^{(n-1)}+\cdots+a_{n-1}y'+a_ny=0$.
\par
The functions $s_{ij}(t)$, $j=0,1,\ldots,k_i-1$, $i=1,\ldots,p$,
form a fundamental system of solutions of the equation
$s^{(n)}(t)+a_1s^{(n-1)}(t)+\cdots+a_{n-1}s'(t)+a_ns(t)=0$.
If we fix
$s_0,s_1,\ldots,s_{n-1}$, there exists a unique solution
$s(t)=s_0+s_1t/1!+\cdots+s_{n-1}t^{n-1}/(n-1)!+t^nu(t)$ for some
$u(t)\in\mathbb{C}[[t]]$ (because $s_j=s^{(j)}(0)$) and a
unique system of constants $\gamma_{ij}\in\mathbb{C}$ such that
$c(t)=\sum\gamma_{ij}c_{ij}(t)$.
Hence the partial sum of the
first $n$ summands of the series
\[
(\gamma_{10}+\gamma_{11}t+\cdots+\gamma_{1,k_1-1}t^{k_1-1})
\text{\rm exp}(\lambda_1t)+
\cdots
\]
\[
+(\gamma_{p0}+\gamma_{p1}t+\cdots+\gamma_{p,k_p-1}t^{k_p-1})
\text{\rm exp}(\lambda_pt)
\]
coincides with $s_0+s_1t/1!+\cdots+s_{n-1}t^{n-1}/(n-1)!$. If we
compare
the coefficients of $t^j$, $j=0,1,\ldots,n-1$,
we obtain a system of $n$ linear equations with unknowns $\gamma_{ij}$,
and the determinant of the system is different from 0.
\par
By Proposition \ref{linear equations}, the only solution
$y(x)=\sum_{j\geq 0}c_j(x)\rho^j$ with
$c_0(x)=c_1(x)=\cdots=c_{n-1}(x)=0$ is the zero formal power
series. Hence, in order to see that an arbitrary solution is a
linear combination of solutions of the form $y_{ij}(c_{ij},x)$, it
is sufficient to show that for any $n$ constants
$c_0(x),c_1(x),\ldots,c_{n-1}(x)\in\widehat{R_0}$ there exist
constants $c_{ij}(x)$ such that the partial sum of the first $n$
summands of the series
\[
(c_{10}+c_{11}\rho+\cdots+c_{1,k_1-1}\rho^{k_1-1})
\text{\rm exp}(\lambda_1\rho)+
\cdots
\]
\[
+(c_{p0}+c_{p1}\rho+\cdots+c_{p,k_p-1}\rho^{k_p-1})
\text{\rm exp}(\lambda_p\rho)
\]
coincide with $c_0+c_1\rho+\cdots+c_{n-1}\rho^{n-1}$. Again, we
consider
the coefficients from $\widehat{R_0}$ of $\rho^j$, $j=0,1,\ldots,n-1$,
form a system of $n$ linear equations with unknowns $c_{ij}$, and
the determinant of this system is the same as the determinant of the
corresponding system for $\mathbb{C}[[t]]$.
Since this determinant is nonzero, the system has a unique solution
and we can find the desired constants $c_{ij}(x)\in \widehat{R_0}$.
\end{proof}

A special case of our considerations is the non-associative
exponential function introduced in \cite{DG}. We define the
exponent $E(x)$ as the formal power series in
$K\{\{x\}\}=\widehat{K\{x\}}$ satisfying the conditions
$E'(x)=E(x)$, $E(0)=1$ and $E(x)E(x)=E(2x)$. By Theorem
\ref{homogeneous case}, all solutions of the equation $E'(x)=E(x)$
are $c(x)\text{\rm exp}(\rho)$, where
$c(x)\in\widehat{(K\{x\})_0}$. There are many solutions satisfying
the condition $E(0)=1$ (which simply means that $c(0)=1$), and
only the second condition $E(x)E(x)=E(2x)$ determines $E(x)$ in a
unique way.

\section*{Acknowledgements}

The first author wants to thank for the hospitality at the
Department of Mathematics of Ruhr University in Bochum during his
visit there, and the second author wants to thank for the
hospitality of the Institute of Mathematics and Informatics of the
Bulgarian Academy of Sciences. Both authors thank L.\ Gerritzen
for stimulating discussions about the subject.

\end{document}